\expandafter\ifx\csname amssym.def\endcsname\relax \else\endinput\fi
\expandafter\edef\csname amssym.def\endcsname{%
       \catcode`\noexpand\@=\the\catcode`\@\space}
\catcode`\@=11

\def\undefine#1{\let#1\undefined}
\def\newsymbol#1#2#3#4#5{\let\next@\relax
 \ifnum#2=\@ne\let\next@\msafam@\else
 \ifnum#2=\tw@\let\next@\msbfam@\fi\fi
 \mathchardef#1="#3\next@#4#5}
\def\mathhexbox@#1#2#3{\relax
 \ifmmode\mathpalette{}{\m@th\mathchar"#1#2#3}%
 \else\leavevmode\hbox{$\m@th\mathchar"#1#2#3$}\fi}
\def\hexnumber@#1{\ifcase#1 0\or 1\or 2\or 3\or 4\or 5\or 6\or 7\or 8\or
 9\or A\or B\or C\or D\or E\or F\fi}

\font\tenmsa=msam10
\font\sevenmsa=msam7
\font\fivemsa=msam5
\newfam\msafam
\textfont\msafam=\tenmsa
\scriptfont\msafam=\sevenmsa
\scriptscriptfont\msafam=\fivemsa
\edef\msafam@{\hexnumber@\msafam}
\mathchardef\dabar@"0\msafam@39
\def\dashrightarrow{\mathrel{\dabar@\dabar@\mathchar"0\msafam@4B}}
\def\dashleftarrow{\mathrel{\mathchar"0\msafam@4C\dabar@\dabar@}}

\def\ulcorner{\delimiter"4\msafam@70\msafam@70 }
\def\urcorner{\delimiter"5\msafam@71\msafam@71 }
\def\llcorner{\delimiter"4\msafam@78\msafam@78 }
\def\lrcorner{\delimiter"5\msafam@79\msafam@79 }
\def\yen{{\mathhexbox@\msafam@55 }}
\def\checkmark{{\mathhexbox@\msafam@58 }}
\def\circledR{{\mathhexbox@\msafam@72 }}
\def\maltese{{\mathhexbox@\msafam@7A }}

\font\tenmsb=msbm10
\font\sevenmsb=msbm7
\font\fivemsb=msbm5
\newfam\msbfam
\textfont\msbfam=\tenmsb
\scriptfont\msbfam=\sevenmsb
\scriptscriptfont\msbfam=\fivemsb
\edef\msbfam@{\hexnumber@\msbfam}
\def\Bbb#1{{\fam\msbfam\relax#1}}
\def\widehat#1{\setbox\z@\hbox{$\m@th#1$}%
 \ifdim\wd\z@>\tw@ em\mathaccent"0\msbfam@5B{#1}%
 \else\mathaccent"0362{#1}\fi}
\def\widetilde#1{\setbox\z@\hbox{$\m@th#1$}%
 \ifdim\wd\z@>\tw@ em\mathaccent"0\msbfam@5D{#1}%
 \else\mathaccent"0365{#1}\fi}
\font\teneufm=eufm10 scaled 1200
\font\seveneufm=eufm7
\font\fiveeufm=eufm5
\newfam\eufmfam
\textfont\eufmfam=\teneufm
\scriptfont\eufmfam=\seveneufm
\scriptscriptfont\eufmfam=\fiveeufm

\csname amssym.def\endcsname

\newcommand{\ext}[1]{\stackrel{#1}{\wedge}}
\def\R{{\Bbb R}}

\def\Z{{\Bbb Z}}
\def\Q{{\Bbb Q}}

\def\X{{{\tilde X}}}

\def\tg{{{\tilde g}}}

\def\l{\lambda }
\def\La{\Lambda }
\def\a1{${\alpha}_1$ }
\def\a2{${\alpha}_2$ }
\def\a{\alpha }

\def\ga{\gamma }
\def\Ga{\Gamma }

\def\r{\rho }
\def\p{\phi }
\def\tp{{\tilde \phi }}

\def\om{\omega }
\def\tom{{\tilde \omega }}
\def\M{{\cal M}_g}
\def\Mg{{\cal M}_g}
\def\Mga{{\cal M}_{\ga}}
\def\Mgk{{\cal M}_g[k]}
\def\Mgr{{\cal M}_g^r}
\def\Mgrk{{\cal M}_g^r[k]}
\def\Mt{{\cal M}_{\tilde g}}
\def\N{{\cal M}_{\tilde g}}

\def\O{{\cal O}}
\def\Pa{Pic_{alg}(\M)}
\def\Paq{{Pic_{alg}(\M)} \otimes {\Q}}
\def\Ph{Pic_{hol}(\M)}
\def\Phq{{Pic_{hol}(\M)} \otimes {\Q}}
\def\Pf{Pic_{fun}(\M)}
\def\Pfr{Pic_{fun}({\M}^{r})}

\def\Pt{Pic(\Tg)}
\def\Pth{Pich(\Tg)}
\def\Ptt{Pich(\Tt)}
\def\Pthq{{Pich(\Tg)}_{\Q}}
\def\Pthqr{{Pich({\Tg}^{r})}_{\Q}}
\def\Pttq{{Pich(\Tt)}_{\Q}}
\def\Pttqr{{Pich({\Tt}^{r})}_{\Q}}
\def\Tg{{\cal T}_g}
\def\Tgr{{\cal T}_g^r}
\def\Tt{{\cal T}_{\tilde g}}
\def\Tga{{\cal T}_{\ga}}
\def\Tin{{\cal T}_{\infty}}
\def\Tir{{\cal T}_{\infty}^{r}}

\def\mapright#1{\smash{
    \mathop{\longrightarrow}\limits^{#1}}}
\def\mapdown#1{\Big\downarrow
 \rlap{$\vcenter{\hbox{$\scriptstyle#1$}}$}}
\documentstyle[12pt]{article}
\begin{document}
\baselineskip17pt
\begin{center}
{\bf{UNIVERSAL HODGE BUNDLE AND MUMFORD ISOMORPHISMS ON THE
ABELIAN INDUCTIVE LIMIT OF TEICHM\"ULLER SPACES}} \\
\vspace{0.4cm}
{\bf Indranil Biswas, Subhashis Nag and Dennis Sullivan} \\
\end{center}
\vspace{0.4cm}
\begin{center}
{\bf Abstract}
\end{center}
\vspace{0.1cm}

Let $\Mg$ denote the moduli space of compact Riemann surfaces
of genus $g$.  Mumford had proved, for each fixed genus $g$, that there are
isomorphisms asserting that certain higher $DET$ bundles over
$\Mg$ are certain fixed (genus-independent) tensor powers of the Hodge line
bundle on $\Mg$. We obtain a coherent, genus-independent
description of the Mumford isomorphisms over certain infinite-dimensional
``universal'' parameter spaces of compact Riemann
surfaces (with or without marked points).
We work with an inductive limit of Teichm\"uller
spaces comprising complex structures on a certain ``solenoidal
Riemann surface'', $H_{\infty,ab}$, which appears as the inverse limit
of an inverse system of surfaces of different genera
connected by abelian covering maps. We construct the universal Hodge and
higher $DET$ line bundles on this  direct limit of
Teichm\"uller spaces (in the sense of Shafarevich). The main result shows
how such $DET$ line bundles on the direct limit carry coherently-glued
Quillen metrics and are related by the appropriate Mumford isomorphisms.

Our work can be viewed as a contribution to a non-perturbative
formulation of the Polyakov measure structure in a genus-independent fashion.

\vskip 0.1cm
\noindent
\vskip 0.1cm
{\bf Authors' addresses:}
\noindent
\vskip 0.1cm
ICTP, P.O. Box 586, 34100-Trieste, ITALY;
\noindent
\vskip 0.1cm
The Institute of Mathematical Sciences, Madras-600 113, INDIA; Max-Planck
Institut Math, 53225 Bonn, GERMANY.
\noindent
\vskip 0.1cm
IHES, Bures-sur-Yvette, 91440-FRANCE; Einstein Chair, CUNY,
New York-10036, USA.

\newpage
\noindent{\underline {\it INTRODUCTION:}}

\noindent
The chief purpose of the present paper is to obtain a genus-independent
description of the Mumford isomorphisms over the moduli spaces of Riemann
surfaces (with or without marked points). Mumford proved that,
for each fixed genus $g$, there are
isomorphisms asserting that certain higher $DET$ bundles over the moduli space
$\M$ are certain fixed (genus-independent) tensor powers of the Hodge line
bundle on $\M$. (Hodge itself is the first of these $DET$ bundles, and
generates the Picard group of $\M$.)

In order to establish a coherent genus-independent description of the
Mumford isomorphisms, we work with an inductive limit of Teichm\"uller
spaces corresponding to a projective family of surfaces of different genera
connected by abelian covering maps. Since a covering space between surfaces
induces (by a contravariant functor) an immersion
of the Teichm\"uller space of
the surface of lower genus into the Teichm\"uller space of the covering
surface, we can think of line bundles living on the direct limit of these
Teichm\"uller spaces (in the sense of [Sha]). The main result shows
how such ``$DET$'' line bundles on the direct limit carry coherently-glued
Quillen metrics and are related by the appropriate Mumford isomorphisms;
these isomorphisms actually become unitary isometries
with respect to the Quillen hermitian structures.

All this is, of course, closely related to the Polyakov bosonic string
theory, -- and our work can be viewed as an attempt at a non-perturbative
formulation of the Polyakov measure structure in a genus-independent mode.
Regarding the Polyakov measure (which is a measure on each moduli space
$\M$ and describable via the Mumford isomorphisms), and its relation to our
work, we shall have more to say in Section 5c.

\noindent
{\bf Acknowledgements:} We would sincerely like to
thank E.Arbarello, E.Looijenga and
M.S.Narasimhan for their interest and many helpful discussions,
and D.Prasad for pointing out to us the paper of Moore [M].
We specially thank Professor Looijenga for generously supplying us,
(in early January 1994),
with a proof of a strong form of the discreteness
of the Picard group of moduli spaces of surfaces with marked points and
level-structure. That result is needed critically in our construction, and
we have inserted the proof he sent us as an Appendix to this work.
Later on he also pointed out to us that the discreteness result
can be found in a recent preprint of R.Hain (see References).

The first two authors would like to thank RIMS/Kyoto University, TIFR
Bombay, and ICTP Trieste, (in chronological order), for creating
very pleasant opportunities for us
to get together and collaborate on this project.

\bigskip
\noindent {\underline {\it Section 1. BASIC BACKGROUND:}}

\noindent {\bf 1a. $\Tg$ and $\Mg$:}
Let $X$ be a compact oriented $C^{\infty}$ surface of genus $g \geq 2$.
Let $T(X)=\Tg$ denote the corresponding Teichm\"uller space and $\M$
denote the moduli space of all complex structures on $X$. $\M$ is a
complex analytic V-manifold of dimension $3g-3$ that is known to be a
quasi-projective variety. $\Tg$ itself can be constructed as the orbifold
universal covering of $\Mg$ in the sense of Thurston; conversely, $\Mg$ can be
recovered from $\Tg$ as the quotient by the action of the
mapping class group (=``Teichm\"uller modular group''), $MC_g$, acting
biholomorphically and properly discontinuously on the contractible, Stein,
bounded, $(3g-3)$-dimensional domain, $\Tg$. In the sequel, unless
otherwise stated, we will assume $g \geq 3$.

\noindent
{\it Remark:} For simplicity, in the main body of the paper we
deal with the situation of closed Riemann surfaces without punctures
and prove our main theorems in this context. In Section 6
we shall show how to extend our main results also to the case of moduli
of Riemann surface with some fixed number of marked points.

\noindent {\bf 1b. Bundles on moduli:}
Our first critical need is to study and compare the Picard groups of
algebraic/holomorphic line bundles over $\Mg$. We introduce three closely
related groups of bundles:

$$
\Pa
= the~group~of~algebraic~line~bundles~over~\M.
\leqno(1.1)
$$

$$
\Ph = H^{1}(\M,{\O}^{*})
=the~group~of~holomorphic~line~bundles~over~\M
\leqno(1.2)
$$

$$
\Pf = the~Picard~group~of~the~moduli~functor.
\leqno(1.3)
$$

\noindent
An element of $\Pf$ consists of the prescription of a holomorphic
line bundle $L_F$ on the base space $S$ for every complex analytic
family of Riemann surfaces of genus $g$ ($C^{\infty}$ locally trivial
Kodaira-Spencer family) $F=(\gamma: V \rightarrow S)$ over any complex
analytic base $S$. Moreover, for every commutative diagram of
families $F_{1}$ and $F_{2}$ having the morphism $\a$ from the base
$S_{1}$ to $S_{2}$, there must be assigned a corresponding isomorphism
between the line bundle $L_{F_{1}}$ and the pullback via $\a$ of the bundle
$L_{F_{2}}$. For compositions of such pullbacks these isomorphisms between
the prescribed bundles must satisfy the obvious compatibility condition.
See [Mum], [HM], [AC].
[Note: Mumford has considered this Picard group of the moduli
functor also over the Deligne-Mumford compactification of $\M$.]

There is, of course, a natural map of $\Pa$  to $\Ph$ by stripping off the
algebraic structure and remembering only the holomorphic structure. Also it
is clear that any actual line bundle over $\M$ produces a corresponding
element of $\Pf$ simply by pulling back to the arbitrary base $S$ using the
canonical classifying map for the given family $F$ over $S$. In fact,
Mumford has shown that $\Pa$ is a subgroup of finite index in $\Pf$, and
that the latter is torsion-free. Therefore $\Pa$ is also torsion free.
See [Mum] and [AC].

Now, it is known that the cohomology $H^{1}(\M, \O)$ vanishes.
[In fact, by Kodaira-Spencer theory of infinitesimal deformations of
structure, we know that the above cohomology represents
the tangent space to $\Ph$. Since by Mumford and Harer's results
(see below) we know that Pic is discrete, this cohomology must vanish.]
It therefore follows that the Chern class homomorphism
$c_{1}: \Ph \rightarrow H^{2}(\M,\Z)$ is an
injection. Since Harer ([Har]) has shown that the $H^{2}(MC_{g},\Z)$
is $\Z$, it follows that $\Ph$ is also torsion free, and in fact therefore
$\Ph$ must be isomorphic to $\Z$ (for $g \geq 5$).

The torsion-free property guarantees that we stand to lose nothing by
tensoring with $\Q$; namely, the maps:
$$
\Pa \rightarrow \Paq ~~~~~~~~~ \Ph \rightarrow \Phq
\leqno(1.4)
$$
are {\it injective}. (All tensor products are over $\Z$.)
We can thus think of the original Picard groups
respectively as {\it discrete lattices} embedded inside these $\Q$-vector
spaces. These two {\it rational Picard groups} above are of
fundamental interest for us. We wish to show that they are naturally
isomorphic:

\noindent {\bf Lemma 1.1:} {\it The structure-stripping map (that forgets the
algebraic nature of the transition functions for bundles in $\Pa$), induces a
canonical isomorphism $\Paq \longrightarrow \Phq$.}

\noindent
{\bf Proof:} The Chern class homomorphism provides us with two maps:
$$
c_a: \Paq \longrightarrow H^{2}(\M,\Q)
\leqno(1.5)
$$

$$
c_h: \Phq \longrightarrow H^{2}(\M,\Q)
\leqno(1.6)
$$
which, by naturality, give a commutative triangle when combined
with the structure-stripping map of the Lemma. The Lemma will therefore
be proved if we show that each of the two Chern class maps above are
themselves isomorphisms.

Since, as we mentioned above, $H^{1}(\M,\O)=\{0\}$, we deduce that
both $c_a$ and $c_h$ are injections.
To prove that they are surjective it clearly suffices to find {\it some}
algebraic line bundle over $\M$ with a non-vanishing Chern class
(because, by Harer's result we know that $H^{2}(\M,\Q) = \Q$.)

Towards this aim, and since it is of basic importance to us, we introduce
the {\it Hodge line bundle, $\l$}, as a member of $\Pf$. This is easy.
Consider a family of genus $g$ Riemann surfaces,
$F=(\gamma: E \rightarrow S)$ over $S$. Then the
``Hodge bundle'' on the parameter space $S$ is defined to be
$det(R^1{\ga}_*{\O} -R^0{\ga}_*{\O}))$. Here the $R^i$ denote the usual
direct image and higher derived functors (see, for example, [H]).
One thus obtains an element of $\Pf$ by definition. It is a known fact that
$\Pf$ is {\it generated} by $\l$; moreover, for $g\geq 3$, $\Pf=\Z$.
(See Theorem 1 [AC].)

Now recall Mumford's assertion that $\Pa$ is a subgroup of finite
index within $\Pf=H^{2}(MC_{g}, \Z)$ ([Mum]). Therefore some integer
multiple (i.e., tensor power) of $\l$ must lie in $\Pa$ and
have non-zero Chern class. We are therefore through. $\hfill{\Box}$

\noindent {\it Remark:}~ By Mumford, some tensor power of each element
of $\Pf$ is an honest algebraic line bundle on moduli.
As a corollary of our Lemma we see that {\it some power every holomorphic
line bundle over moduli space carries an algebraic structure.}

\smallskip
\noindent
{\it Discussion of the Hodge line bundle:}
Let ${\M}^{\!\!o}$ denote the open set of moduli space comprising Riemann
surfaces without nontrivial automorphisms.
Let us utilise the universal genus $g$ family over $\Tg$,
and restrict the family over ${\Tg}^{o}$ (automorphism-free points of
Teichm\"uller space). Clearly then, by quotienting out
the fix-point free action of $MC_{g}$ on this open subset of $\Tg$, we
obtain over ${\M}^{\!\!o}$
a Hodge line bundle whose fiber over the Riemann surface $X$ is
$\ext{g}H^0(X,K_X)^*$. (We let $K_X$ denote the canonical
line bundle of $X$.) What we see from the proof is that although the Hodge
line bundle may not exist on the whole of $\M$ as an actual (algebraic or
holomorphic) line bundle, some integer power of the Hodge bundle over
${\M}^{\!\!o}$ necessarily lives as an actual algebraic line bundle over the
full moduli space. {\it Therefore, when we go to the rational Picard groups
over $\M$, the Hodge itself exists as an element of $\Paq = \Phq$.} This is
of crucial importance for our work.

\medskip
Our chief object in making the above careful considerations is to be able
to think of the Picard group over moduli as a group of {\it
$MC_g$-invariant line bundles on the Teichm\"uller space}. Indeed,
there is a 1-1 correspondence between the Picard group $\Ph$
on $\Mg$ and the set of pairs of the form $(L,\psi)$, where
$L$ is a holomorphic line bundle on $\Tg$ and $\psi$ is a lift of the
action $MC_g$ (the mapping class group) on $\Tg$ to the the
total space of $L$. The pull-back of a bundle on $\Mg$ to $\Tg$ is
naturally equipped with a lift of the action of $MC_g$. Conversely, a
line-bundle on $\Tg$ equipped with a lift of the action of $MC_g$
descends to a line bundle on $\Mg$. Define

$$
\Pt\,:=\,\ \{(L,\psi)\,|~~ L ~{\rm and}~ \psi ~{\rm as~ above}\}
\leqno(1.7)
$$

We will refer to the elements of $Pic(\Tg)$ as {\it
modular invariant line-bundles}. Two elements $(L,\psi)$ and
$(L',{\psi}')$ are called isomorphic if there is an isomorphism $h:\,
L\rightarrow\,L'$ which is equivariant with respect to the action of
$MC_g$. In other words, ${\psi}'\circ h=h\circ\psi$.

$Pic(\Tg)$ has a natural structure of an abelian group, and, from the
above remarks, is isomorphic as a group to $\Ph$.

Now we introduce hermitian structure on elements of $Pic(\Tg)$. Let
$\Pth$ be the set of triplets of the form $(L,\psi,h)$, where
$(L,\psi) \in Pic(\Tg)$ and $h$ is a hermitian metric on $L$ which is
invariant under the action $MC_g$ given by $\psi$. Note that a
modular invariant metric $h$, as above, induces a modular invariant
metric on any power of $L$. Hence $\Pth$ has a structure of an
abelian group.

\noindent
{\bf 1c. DET bundles for families of d-bar operators:} Given, as before,
any Kodaira-Spencer family $F=(\gamma: V \rightarrow S)$,
of compact Riemann surfaces of genus $g$, and a holomorphic vector
bundle $E$ over the total space $V$, we can consider the base $S$ as
parametrizing a family of elliptic d-bar operators, as is standard.
The operator corresponding to $s \in S$ acts along the fiber Riemann
surface $X_s = {\gamma}^{-1}(s)$ :
$$
\bar{\partial}_{s}: C^{\infty}({\gamma}^{-1}(s), E) \rightarrow
C^{\infty}({\gamma}^{-1}(s), E \otimes {{\Omega}^{0,1}_{X_{s}}})
\leqno(1.8)
$$
One defines the associated vector space of one dimension given by:
$$
DET(\bar{\partial}_{s}) = (\ext{max}ker \bar{\partial}_{s})^{*} \otimes
(\ext{max} coker \bar{\partial}_{s})
\leqno(1.9)
$$
and it is known that these complex lines fit together naturally over the
base space $S$ giving rise to a holomorphic line bundle over $S$ called
$DET(\bar{\partial})$. In fact, this entire construction is natural with
respect to morphisms of families and pullbacks of vector bundles.

We could have followed the above construction through for the universal
genus $g$ family over $\Tg$, with the vector bundle $E$ being, for example,
the trivial line bundle over the universal curve, or
the vertical (relative) tangent bundle, or any of its higher tensor powers.
It is clear that setting $E$ to be the trivial line bundle over
$V$ for any family $F=(\gamma: V \rightarrow S)$, the above prescription
for $DET$ provides merely another
description of the Hodge line bundle over the base $S$. Indeed, the fiber
of $DET=Hodge$ over $s$ is $\ext{g}H^0(X_{s},K_{X_{s}})^*$, and the
naturality of the $DET$ construction mentioned above shows that we thus get a
member of $\Pf$.

By the same token, setting over any family $F$ the vector bundle $E$ to be
the $m^{th}$ tensor power of the vertical tangent bundle along the fibers,
we get by the DET construction a well-defined member
$$
{\l}_{m} \in \Pf, ~~~ m=0,1,2,3,...
\leqno(1.10)
$$
\noindent
Clearly, $\l_{0}$ is the Hodge bundle in the Picard group of the moduli
functor. By ``Teichm\"uller's lemma'' one notes that $\l_{1}$ represents
the canonical bundle over the moduli space, the fiber at any Riemann
surface $X$ being the determinant line of the space of holomorphic
quadratic differentials thereon.

\noindent
{\bf 1d. Mumford isomorphisms:} By applying the Grothendieck-Riemann-Roch
theorem it was proved ([Mum]) that as elements of $\Pf$ one has the following
isomorphisms:
$$
{\l}_m = (6m^{2} + 6m + 1)~~ tensor~~power~~of~~Hodge(\l_{0})
\leqno(1.11)
$$

It is well-known, (Satake compactification combined with Hartogs theorem),
that there are no non-constant holomorphic functions on $\M$ ($g \geq 3$).
Therefore the choice of an isomorphism of ${\l}_{m}$ with
${{\l}_{0}}^{\otimes {(6m^{2}+6m+1)}}$ is {\it unique} up to a nonzero
scalar. We would like to put canonical Hermitian metrics on these DET
bundles so that this essentially unique isomorphism actually becomes an
unitary isometry. This is the theory of the:

\noindent
{\bf 1e. Quillen metrics on DET bundles:}
If we prescribe a conformal Riemannian metric on the fiber Riemann
surface $X_s$, and simultaneously a hermitian fiber metric on the vector
bundle $E_s$, then clearly this will induce a natural $L^2$ pairing on the
one dimensional space $DET(\bar{\partial}_{s})$ described in (1.9). Even if
one takes a smoothly varying family of conformal Riemannian metrics on the
fibers of the family, and a smooth hermitian metric on the vector bundle
$E$ over $V$, these $L^2$ norms on the DET-lines may fail to fit together
smoothly (basically because the dimensions of the kernel or cokernel for
${{\bar{\partial}}_{s}}$ can jump as $s$ varies over $S$). However,
Quillen, and later Bismut-Freed and other authors, have described
a ``Quillen modification'' of the $L^2$ pairing which always produces a
smooth Hermitian metric on $DET$ over $S$, and has functorial properties.
[Actually, in the cases of our interest the Riemann-Roch theorem
shows that the dimensions of the kernel and cokernel spaces
remain constant over moduli -- so that the $L^2$ metric is itself
smooth. Nevertheless the Quillen metrics will be crucially used by us
because of certain functorial properties, and curvature properties,
that they enjoy.]

In fact, using the metrics assigned on the Riemann surfaces (the fibers of
$\gamma$), and the metric on $E$, one gets $L^2$ structure on the spaces of
$C^{\infty}$ sections that constitute the domain and target for our d-bar
operators. Hence ${\bar{\partial}}_{s}$ is provided with an adjoint
operator ${{\bar{\partial}}_{s}}^{*}$, and one can therefore construct the
positive (Laplacian) elliptic operator as the composition:
$$
{\Delta}_{s} =  {{\bar{\partial}}_{s}}^{*} \circ {{\bar{\partial}}_{s}},
$$
These Laplacians have a well-defined (zeta-function regularized)
determinant, and one sets:
$$
Quillen~norm~on~fiber~of~DET~above~s
= (L^{2}~norm~on~that~fiber)(det{\Delta}_{s})^{1/2}
\leqno(1.12)
$$
{\it This turns out to be a smooth metric on the line bundle $DET$.}
See [Q], [BF].

In the situation of our interest, the vector bundle $E$ is the vertical
tangent (or cotangent) line bundle along the fibers of $\gamma$, or its
powers, so that the assignment of a metric on the Riemann surfaces already
suffices to induce a hermitian metric on $E$. Hence one gets a Quillen norm
on the various $DET$ bundles ${\l}_m$ ($\in \Pf$) for every choice
of a smooth family of conformal metrics on the Riemann surfaces.
{\it The Mumford isomorphisms (over any base $S$) become isometric
isomorphisms with respect to the Quillen metrics.}

The {\it curvature form} (i.e., first Chern form) on the base $S$
of the Quillen DET bundles has a particularly elegant expression:
$$
c_{1}(DET, Quillen~metric) = - \int_{V|S}^{} (Ch(E) Todd({T}_{vert}))
\leqno(1.13)
$$
\noindent
where the integration represents integration of differential forms along
the fibers of the family $\gamma: V \rightarrow S$.

We now come to one of our main tools in this paper.
By utilising the uniformisation theorem (with moduli parameters), the
universal family of Riemann surfaces over $\Tg$, and hence any holomorphic
family $F$ as above, has a smoothly varying family of Riemannian metrics on
the fibers given by the {\it constant curvature -1 Poincare
metrics}. {\it The
Quillen metrics arising on the $DET$ bundles $\l_m$ from the Poincare
metrics on $X_s$ has the following beautiful property for its curvature}:
$$
c_1({\l}_{m}, Quillen) = (12{\pi}^2)^{-1}(6m^{2}+6m+1){{\omega}_{WP}},
{}~~m=0,1,2,..
\leqno(1.14)
$$
\noindent
where ${\omega}_{WP}$ denotes the (1,1) K\"ahler form on $\Tg$ for the
classical {\it Weil-Petersson} metric of $\Tg$. (We remind the reader that
the cotangent space to the Teichm\"uller space at $X$ can be canonically
identified with the vector space of holomorphic quadratic differentials
on $X$, and the WP Hermitian pairing is obtained as
$$
(\phi, \psi)_{WP} = \int_{X}^{} \phi {\bar {\psi}} (Poin)^{-1}
\leqno(1.15)
$$ (Here $(Poin)$ is the area
form on $X$ induced by the Poincare metric.) That the curvature formula
(1.13) takes the special form (1.14) for the Poincare family of metrics has
been shown by Zograf and Takhtadzhyan [ZT]. Indeed, (1.13) specialised to
$E=T_{vert}^{\otimes m}$ becomes simply $-(6m^{2}+6m+1)/12$ times
${\int_{V|S}^{}}{c_{1}{(T_{vert})}^{2}}$. This last integral
represents, for the Poincare-metrics family, ${\pi}^{-2}$ times
the Weil-Petersson symplectic form. See also [BF], [BK], [Bos].

Applying the above machinery, we will investigate in this paper the
behaviour of the Mumford isomorphisms in the situation of a covering map
between surfaces of different genera. Since the Weil-Petersson form is {\it
natural} with respect to coverings, the above facts will be very useful.


\bigskip
\noindent {\underline{\it Section 2. COVERINGS AND DET BUNDLES:}}

\noindent
{\bf 2a. Induced morphisms on moduli:}
Over the genus $g$ closed oriented surface $X$ let:
$$
\pi: \X \longrightarrow\, X
\leqno(2.1)
$$
be  any unramified covering space, where $\X$ is a surface of genus $\tg$.

If $X$ is equipped with a complex structure then there is exactly one
complex structure on $\X$ such that the map $\pi$ is a holomorphic
map. In other words, a complex structure on $X$ induces a complex
structure on $\X$. One may tend to think that this correspondence
induces a morphism from $\Mg$ to $\Mt$. But a closer analysis shows
that this is not in general the case. In fact, the space $\Mg$ can be naturally
identified with $Comp/Diff^{+}$, where $Comp$ is the space of all complex
(=conformal) structures on $X$ and $Diff^{+}$ is the group of all
orientation-preserving diffeomorphisms of $X$.  But
arbitrary element of $Diff^{+}$ can not be lifted to a diffeomorphism of
$\X$ -- and this is the genesis of why there is in general no induced map
of $\Mg$ to $\Mt$.

Since the Teichm\"uller space $\Tg$ is canonically identified with the
quotient $Comp/Diff_{0}$, (where $Diff_{0} \subset Diff^{+}$ is
the sub-group consisting of all those diffeomorphisms which are homotopic
to the identity map), we see however that there is in fact a map induced by
the covering between the relevant Teichm\"uller spaces. Indeed,
using the homotopy lifting property, a diffeomorphism  homotopic to
identity lifts to a diffeomorphism of $\X$; moreover, the lifted
diffeomorphism is also homotopic to the identity. Therefore,
{\it any covering $\pi$ induces a natural morphism between the complex
manifolds $\Tg$ and $\Tt$ obtained by ``pullback'' of complex structure:}
$ Teich(\pi) = {\pi}_T: \Tg \longrightarrow \Tt $.
It is well-known that the induced map is a holomorphic injective immersion
that preserves the Teichm\"uller distance.

However, for our work we need to compare certain line bundles at
an appropriate moduli-space (or finite covering thereof) level, where we
can be sure that the space of line bundles is {\it discrete}. Therefore it
is imperative that we find whether the covering $\pi$ induces a map at the
level of some finite coverings of moduli space. To investigate this -- and
give an affirmative answer in certain topologically restricted situations
-- we introduce the finite coverings of moduli space called:

\noindent
{\it Moduli of Riemann surfaces with level structure:}
Fix any integer $k$ at least two. Let $\Mgk$ denote the moduli with
level $k$. In other words, $\Mgk$ is the parameter space of pairs of the form
$(M,\r)$, where $M$ is a Riemann surface of genus $g$,
and $\r$ is a basis of the $\Z/k\Z$ module $H_1(M,\Z/k\Z)$.

Now, Teichm\"uller space parametrises pairs $(M, \r)$
where $\r$ is a presentation of the fundamental group of the Riemann
surface $M$. Hence it is clear that $\Mgk$ is
a quotient of the Teichm\"uller space by a normal subgroup of the Teichm\"uller
modular group. $\Mgk$ is a finite branched cover over $\Mg$,
$\psi: \Mgk \rightarrow \Mg$. The covering group is $Sp(2g, \Z/k)$.

Now suppose that the covering space $\pi$ is a {\it cyclic} covering:
namely, $\pi\,:\,\ \X\,\longrightarrow\, X$ be a Galois cover, with the
Galois group being $\Z/l\Z$, where $l$ is a {\it prime} number.
Of course, the number of sheets $l$ is the ratio of the Euler
characteristics, i.e., $(\tg - 1)/(g - 1)$.
{\it Fix $k$ any positive integer that is a multiple of the
sheet number $l$.} Under these special circumstances we  show:

\smallskip
\noindent {\bf Proposition\,2.1.}\,\ {\it The covering $\pi$ induces
natural maps
$$(i)\quad\,\ {\pi}_T\,:\,\ \Tg\,\longrightarrow\,\Tt$$

$$(ii)\quad\,\ {\pi}_k\,:\,\ \Mgk\,\longrightarrow\,\Mt$$
Moreover, if
$q:\,\Mgk\longrightarrow\,{\cal M}_g^n[l]$ is the projection induced by
the homomorphism $\Z/k\rightarrow \Z/l$ then the following diagram
commutes $$\matrix{\Mgk&\mapright{id}&\Mgk\cr\mapdown{q} &&
\mapdown{\pi_k}\cr {\cal M}_g[l]&\mapright{\pi_l}&\Mt\cr}$$}
\smallskip

\noindent
{\bf Proof:} Part (i) we have already explained above (without any
topological restrictions on the covering).

In order to prove that $\pi$ induces a map from $\Mgk$ to $\Mt$
note that $\Mgk$ can be identified with $Comp/D(k)$, where $D(k)\subset
Diff^{+}$ denotes the subgroup consisting of all those diffeomorphism whose
induced homomorphism on $H_1(X,\Z/k{\Z})$ is identity. So we need to
check that any element $f\in D(k)$ lifts to a diffeomorphism of $\X$.
{}From the given condition, $\pi$ is an unramified cyclic
cover of order $l$. Now the set of unramified cyclic covers of order
$l$ of $X$ is parametrized by $PH^1(X,{\Z}/l)$. (This denotes the
projective space of the first cohomology vector space over the field
$\Z/l$.) To see this, note that any
such cover $h:Y \rightarrow X$ induces a homomorphism $$h_*\,:\,\
H_1(Y,\Z/l)\,\longrightarrow\, H_1(X,\Z/l)$$ with the quotient being
isomorphic to $\Z/l$. So the quotient map $H_1(X,\Z/l)\rightarrow
H_1(X,\Z/l)/Im(h_*)$ gives an element ${\bar h}\in
PH^1(X,{\Z}/l)$. Conversely, given an element of the $PH^{1}$, there is
a homomorphism $h_*$, and consequently a covering of $X$ determined by
the subgroup of the fundamental group corresponding to $Im(h_{*})$.

Clearly the covering $$f\circ h\, :\,\ Y\,\longrightarrow\,X$$
corresponds to $f^*{\bar h}\in PH^1(X,{\Z}/l)$. But $f^*$ on
$H^1(X,{\Z}/k)$ is the identity, and hence $f^*{\bar h}=\bar h$. We
are through. $\hfill{\Box}$

\noindent {\it Remark:} It is clear that the map $\pi_T$, which exists
for {\it arbitrary} covering $\pi$, will not descend in general to
${\cal M}_g$ because of the non-existence of lifts for arbitrary
diffeomorphisms. That is why we needed to use cyclic coverings and fix a
level structure on the genus $g$ moduli. On the other hand it is worth
noting that the map $\pi_k$ does not respect level structure on the
target side. Namely, even when the covering $\pi$ is cyclic of order $l$
as above, there is no induced map from  $\Mgk$ to ${\Mt}[k]$ (even for
$k=l$). That can be proved by
constructing an example: in fact one can write down a diffeomorphism on
$X$ which acts trivially on homology mod$l$ but such that each lifting to
$\X$ acts nontrivially on the homology mod$l$ up on $\X$.

We also asked ourselves whether the inductive limit space construction,
corresponding to the projective system of coverings
between surfaces, can be carried through at the {\it Torelli spaces}
level. Again it turns out that an arbitrary diffeomorphism that acts
trivially on the homology of $X$ may not lift to such a diffeomorphism on
$\X$.

{\it It therefore appears to be in the very nature of things that we
are forced to do our direct limit construction of spaces and bundles at the
Teichm\"uller level rather than at some intermediate moduli level.}

\noindent
{\bf 2b. Comparison of Hodge bundles:}
We are now in a position to compare the two candidate Hodge bundles
that we get over $\Mgk$; one by the $\pi_k$ pullback of Hodge from $\Mt$,
and the other being its `own' Hodge arising from pulling back Hodge on
$\Mg$.

\noindent {\it Notations:}
{\it Let $\l = {\l}_{0}$ denote, as before, the Hodge bundle on $\Mg$
(a member of $\Pf$, as explained), and let
${\tilde{\lambda}}$ denote the Hodge line
bundle over $\Mt$}.  Also let $\psi: \Mgk \rightarrow \Mg$ denote the
natural projection from moduli with level-structure to moduli.
Let $\om = {\om}_{WP}$ and $\tom$ represent the Weil-Petersson forms (i.e.,
the K\"ahler forms corresponding to the WP hermitian metrics) on
$\Mg$ and $\Mt$, respectively.

The naturality of Weil-Petersson forms under coverings is manifest in the
following basic Lemma:

\smallskip
\noindent {\bf Lemma\, 2.2.}\,\ {\it The 2-forms $l({\psi}^{*}\om)$
and $({\pi}_k)^{*}\tom$ on $\Mgk$ coincide.}
\smallskip

\noindent {\bf Proof.}\,\ This is actually a straightforward
computation. Recall that the cotangent space
to the moduli space is given by the space of quadratic differentials.
Now at any point $\a:=(M,\r)\in \Mgk$, the morphism of cotangent spaces
induced by the map $\pi_k$ of Proposition 2.1(ii) is a map:
$$(d{\pi}_k)^*\,:\,\ T^*_{\pi_k(\a)}\Mt\,\longrightarrow\,T^*_{\a}\Mgk$$
The action on any quadratic differential $\p$ on the Riemann surface
$\pi_k(\a)$, (i.e., $\p \in H^0(\pi_k(\a), K^2)$),  is given by:

$$
(d{\pi}_k)^*{\p}\,=\, 1/l(\sum_{f\in Deck}^{}f^*{\p})
\leqno(2.2)
$$
Here $Deck$ denotes, of course, the group of deck transformations
for the covering $\pi$. Now recall that a covering map $\pi$ induces a
local isometry between the respective Poincare metrics, and that $l$ copies
of $X$ will fit together to constitute $\X$. The lemma therefore follows
by applying formula (2.2) to two quadratic differentials,
and pairing them by Weil-Petersson pairing as per definition (1.15).
$\hfill{\Box}$
\smallskip

The above Lemma, combined with the statements about the curvature form of the
Hodge bundles being $(12{\pi}^{2})^{-1}$ times $\om_{WP}$, show that the
{\it curvature forms (with respect to the Quillen metrics) of the
two line bundles ${\l}^{l}$ and $({\pi}_k)^*{\tilde{\lambda}}$ coincide}.
Do the bundles themselves coincide? Yes:

\smallskip
\noindent {\bf Theorem 2.3.}\,\ {\it The two line-bundles
$({\psi})^*{\l}^l$ and $({\pi}_k)^*{\tilde{\lambda}}$,
on $\Mgk$, are isomorphic.}
{\it If $g\geq 3$, then such an isomorphism:
$$F: ({\psi})^*{\l}^l \longrightarrow ({\pi}_k)^*{\tilde{\lambda}}$$
is uniquely specified up to the choice of a nonzero scaling constant.
The same assertions hold for the higher $DET$ bundles.}
\smallskip

\noindent {\bf Proof.}\,\  Of course, the basic principle of the proof is
that ``curvature of the bundle determines the bundle uniquely'' over the
moduli space $\Mgk$. That would follow automatically if we know that the
Picard group of $\Mgk$ is {\it discrete}. This is in fact known by a recent
preprint of R. Hain, and E. Looijenga has provided us with an
independent proof of (a somewhat stronger) fact. See the Appendix to this
paper.

Without using the Hain-Looijenga proof we can prove what is necessary for
our purposes by using a certain cohomology vanishing theorem for $l$-adic
symplectic groups to be found in work of Moore[M]:

Firstly note that, in view of the commutative diagram in
Proposition 2.1, it is enough to prove things for the case $k=l$. The
Quillen metrics on $\l$ and ${\tilde{\lambda}}$ induce hermitian metrics on
the two bundles under scrutiny: $\xi:=
({\psi})^*{\l}^l$ and ${\xi}':= ({\pi}_{l})^*{\tilde{\lambda}}$.
We know their curvatures coincide.
So the line-bundle ${\xi}^*\otimes \xi'$ with the induced metric is
flat. In other words, after holonomy consideration, ${\xi}^*\otimes\xi'$ is
given by a unitary character $\r$ of the fundamental group
${\Ga}':={\pi}_1({\cal M}_g[l])$.

Let us temporarily denote by
${\Ga}:=MC_g$ the mapping class group for genus $g$. Also let
$G:=Sp(2g,{\Z}/l)$, be the symplectic group defined over the field
${\Z}/l$. We have an exact sequence
$$0\,\longrightarrow\,{\Ga}'\,
\longrightarrow\,\Ga\,\longrightarrow\,G\,\longrightarrow\, 0$$
which induces the following long exact sequence of cohomologies
$$\longrightarrow\,
Hom({\Ga},U(1))\,\longrightarrow\,Hom({\Ga}',U(1))\,\longrightarrow\,
H^2(G,U(1))\,\longrightarrow $$
Note that all the actions of $U(1)$ are taken to be trivial.

Now we claim that in order to prove the theorem it is enough
to show that $H^2(G,U(1))=0$. This is because, in that case the above
exact sequence of cohomologies would imply that $\r$ is the pull-back of
a character on $\Ga$. In other words, ${\xi}^*\otimes\xi'$ would be shown
to be a pull-back of
a flat hermitian line-bundle on $\Mg$. Now, the rational Chern classes of
a flat line bundle vanish, and a line bundle on $\Mg$ with vanishing
1-st Chern class is actually a trivial bundle (Theorem 2, [AC]). Thus
if $H^2(G,U(1))=0$, then the triviality of the line-bundle
${\xi}^*\otimes\xi'$ is established.

So all we have to show is the vanishing of $H^2(G,U(1))$.
That follows from the literature on $l$-adic linear groups. Indeed,
from Chapter III of [M] we get that ${\pi}_1(G)=0$. But then the Lemma
(1.1) of [M] implies that $H^2(G,U(1))=0$.

The uniqueness assertion follows since for $g \geq 3$,
$\Mgk$ does not admit any nonconstant global holomorphic function.
$\hfill{\Box}$
\medskip

\noindent
{\bf 2c. Lifting the identifications to Teichm\"uller space:}
Since the principal level structure keeps changing as we deal with
different covering spaces, we must perforce pullback these bundles to
Teichm\"uller space and use our results above to produce a {\it canonical}
identification of the bundles also at the Teichm\"uller level.

\noindent {\it Notations:}
Let $p:\, \Tg\,\longrightarrow\,\Mg$ denote the quotient projection
from Teichm\"uller to moduli. Similarly, set $\tp:\,\Tt\longrightarrow \N$.
There is the natural map $q:\,\Tg \rightarrow \Mgk$ to the intermediate
moduli space with level-structure, and we have $\psi\circ q=p$.
(Recall, $\psi$ denotes the finite covering from $\Mgk$ onto $\Mg$.)
Furthermore, by functoriality of the operation of pullback of complex
structure, one has ${\pi}_k\circ q=\tp\circ {\pi}_T$.

\noindent {\it Remark:} Since Teichm\"uller spaces are contractible Stein
domains, any two holomorphic line
bundles on $\Tg$ are isomorphic; but, since $\Tg$ admits plenty of
nowhere zero non-constant holomorphic functions, a priori there is
no natural isomorphism. Our problem in manufacturing $DET$ bundles over the
inductive limit of Teichm\"uller spaces is to find natural connecting
bundle maps.

\smallskip
\noindent {\bf Theorem 2.4}
{\it The isomorphism $F$ of Theorem 2.3 induces an isomorphism
of ${\l}^{l}$ and ${\tilde{\lambda}}$ pulled back to $\Tg$:}
$${\bar F}: p^*{\l}^l\,\longrightarrow\, (\pi_T)^* {\tilde{\lambda}}$$
{\it ${\bar F}$ does not depend on the choice of the level $k$.}
\smallskip

\noindent {\bf Proof:} To see how this map $\bar F$ depends upon
$k$ replace $k$ by $k'$, both being multiples of the sheet number $l$.
We may assume $k'$ itself is a multiple of $k$. (If $k'$ is
not a multiple of $k$ we can go to the L.C.M. of the two and
proceed with the following argument.) One has the
following commutative diagram
$$
\matrix{\Tg&\mapright{id}&\Tg\cr
\mapdown{}&&\mapdown{f}\cr {\cal M}_g[k']&\mapright{h}&\Mgk\cr}
$$
where $id$ denotes the identity map and $h$ is the natural
projection. Clearly,  ${\pi}_{k'}= {\pi}_k\circ h$. Therefore
$\bar F$ does not depend upon the choice of $k$.

So we are able to pick up a {\it unique
isomorphism} between $p^*{\l}^l$ and $(\pi_T)^*{\tilde{\lambda}}$, ambiguous
only up to a non-zero scaling constant.
$\hfill{\Box}$

In the next section we will examine some functorial properties of this
chosen isomorphism.


\bigskip
\noindent {\underline{\it Section 3. FROM CYCLIC TO ABELIAN COVERS:}}

We aim to generalise the Theorem 2.4 above to the situation where the
covering $\pi$ admits a factorisation into a finite succession of {\it
abelian} Galois coverings.

\noindent {\bf 3a. Factorisations into cyclic coverings:}
Let $f:\, Z\longrightarrow\,X$ be an {\it abelian} Galois covering
of a surface of genus $g$ by a surface of genus $\ga$, having $N$ sheets.
Suppose that $f$ allows decompositions in two ways as a product of
cyclic coverings of prime order. One therefore has a commutative diagram of
covering maps:

\smallskip
$$\matrix{ Z&\mapright{id} &Z\cr\mapdown{{\r}_1}
&&\mapdown{{\r}_2}\cr Y_1 &&Y_2\cr\mapdown{{\pi}_1}
&&\mapdown{{\pi}_2}\cr X &\mapright{id}&X\cr}
\leqno(3.1)
$$
Here ${\pi}_i \circ {\r}_i =f$, $i=1,2$, and ${\pi}_i$ and
${\r}_i$ are cyclic Galois covers of prime orders $l_i$ and $m_i$
respectively. Let $\a_1$, $\a_2$ be the genera of $Y_1$, $Y_2$
respectively. For the above coverings we have the following
commutative diagram of maps between Teichm\"uller spaces

$$\matrix{\Tg&\mapright{id}&\Tg\cr\mapdown{{\pi}^T_{1}}
&&\mapdown{{\pi}^T_2}\cr
{{\cal T}_{{\a}_1}} && {{\cal T}_{{\a}_2}}
\cr\mapdown{{\r}^T_1}&&\mapdown{{\r}^T_2}\cr \Tga
&\mapright{id}&\Tga\cr}
$$

Now,  $N=l_{1}m_{1} = l_{2}m_{2}$; set
$\pi^T :={\pi}^T_{1}\circ {\r}^T_1={\pi}^T_{2}\circ {\r}^T_2$.

Following the prescription described in the previous section,
corresponding to each of
the two different decomposition of the covering $Z\rightarrow X$, we
will get two isomorphisms of the bundles over $\Tg$:

$$
F_1,~F_2:~\ p^{*}{\l}^N\,\longrightarrow\,{\pi}^{T*}\xi
\leqno(3.2)
$$
where we have denoted by ${\l}$ the Hodge bundle on $\Tg$ and by $\xi$
the Hodge bundle on $\Tga$. For our final construction we require the
crucial Proposition below:

\smallskip
\noindent
{\bf Proposition 3.1.}\,\ {\it The two isomorphisms $F_1$
and $F_2$ are constant (nonzero) multiples of each other.}
\smallskip

\noindent
{\bf Proof:} As in Proposition 2.1, we get the induced maps:
$$
{\pi}_{1,N} ~~({\rm resp.}~\,{\pi}_{2,N}):\,\ {\cal M}_g[N]\,
\longrightarrow \,{{\cal M}_{{\a}_1}} ~~({\rm resp.}~\,{\cal M}_{{\a}_2})
$$
For $i=1,2$, let $q_i$ be the projection of
${\cal M}_{{\a}_i}[m_i]$ onto ${\cal M}_{{\a}_i}$.

The technique of the proof is to go over to bundles on the fiber
product of moduli spaces with level structures and show that isomorphisms
are unique (up to scalar) there. Then the result desired will follow.

The fiber product
${\cal P}:=\,{\cal M}_g[N]{\times}_{({\cal
M}_{{\a}_1}\times {\cal M}_{{\a}_2})}({{\cal
M}_{{\a}_1}}[m_1]\times {{\cal M}_{{\a}_2}}[m_2])$
fits into the following commutative diagram

\smallskip
$$
\matrix{{\cal P}& \mapright{p_1}&{\cal
M}_g[N]\cr\mapdown{p_2}&&\mapdown{({\pi}_{1,l}, {\pi}_{2,l})}\cr
{{\cal M}_{{\a}_1}}[m_1]\times {{\cal
M}_{{\a}_2}}[m_2]&\mapright{(q_1,q_2)}&{\cal M}_{{\a}_1}\times
{\cal M}_{{\a}_2}\cr}
\leqno(3.3)
$$

Again let $f_i:\,\ {{\cal M}_{{\a}_i}}[m_i]\, \longrightarrow\,\Mga$
be the morphisms induced on moduli by the relevant covering spaces,
$(i=1,2)$.  Also let $pr_i$ denote the projection of
${{\cal M}_{{\a}_1}}[m_1] \times {{\cal M}_{{\a}_2}}[m_2]$
to the $i-th$ factor.

In the above notation we have:
$$f_1\circ pr_1\circ p_2\,~ =\,~ f_2\circ pr_2\circ p_2.$$

{}From Theorem 2.3 and diagram (3.3) it now follows that the
following two line bundles over $\cal P$ are isomorphic:
$(f_1\circ pr_1\circ p_2)^*{\xi}$ is isomorphic to
$(\psi\circ p_1)^*{\l}^N$. (Here $\psi$, as in previous sections,
denotes the projection of $\Mg[N]$ to $\Mg$.)

Now, the Teichm\"uller space $\Tg$ has a canonical mapping
to the fiber product $\cal P$.
Clearly the isomorphism $F_i$ under concern ($i=1,2$), is the
pullback of an isomorphism between $(f_i\circ pr_i\circ p_2)^*{ \xi}$
and $(\psi\circ p_1)^*{\l}^N$ over $\cal P$. The desired result will follow
if we show that two such isomorphisms over $\cal P$ can only differ
by a multiplicative scalar.

But that last result is true for any space on which the only
global holomorphic functions are the constants. That is known to be
true for the moduli space $\Mg[N]$. But the fiber product sits over
this by a {\it finite covering} map $p_1$. Therefore the result is
true for the space $\cal P$ as well. (Indeed, if $\cal P$ has a nontrivial
function $f$, then by averaging the powers of $f$ over
the fibers of $p_1$ one would
obtain nontrivial functions on $\Mg[N]$. This argument is well-known.)
The Proposition is proved. $\hfill{\Box}$


\noindent {\bf 3b. Canonical isomorphism from abelian covering:}
The above proposition holds, and the same proof goes through,
if the covering map $f$ is factored by more
than one intermediate covering. In other words, for a situation as below
$$
\matrix{Z&\mapright{}&Y_1^1&\mapright{}&Y_1^2&\ldots
&\mapright{}&X\cr\mapdown{id}&&\mapdown{}&&\mapdown{}&\ldots &&\mapdown{id}\cr
Z&\mapright{}&Y_2^1&\mapright{}&Y_2^2&\ldots &\mapright{}&X\cr}
$$
But the structure theorem for finite abelian group says that any finite
abelian group is a direct sum of cyclic abelian groups whose
orders are prime powers. Moreover any cyclic group of order prime
power $l^n$ fits as the left-end term of an exact sequence of abelian
groups where each quotient group is ${\Z}/l$.
The upshot is that an abelian covering factorises into a sequence
of cyclic covers of prime order.

Utilising Proposition 3.1 and the above remarks we therefore obtain the
sought-for generalisation of Theorem 2.4:

\smallskip
\noindent
{\bf Theorem 3.2:}
{\it Let ${\pi}: Y\,\longrightarrow\,X$ be a covering space, of order
$N$; let $g$ and $\tg$ be the genera of $X$ and $Y$, respectively.
Assume that $\pi$ can be factored
into a sequence of abelian Galois covers. Then there is a canonical
isomorphism (uniquely determined up to a nonzero-constant):

$${\bar F}=\rho: p^*{\l}^N \longrightarrow {\pi}_T^{*}{\tilde{\lambda}}$$

Here ${\pi}_T: \Tg\,\longrightarrow\,\Tt$ is
the map induced by $\pi$ at the Teichm\"uller space level, and ${\l}$ and
$\tilde {\l}$ denote respectively the Hodge bundles on $\Tg$ and $\Tt$.

A similar statement holds using the $m^{th}$ $DET$ bundles instead of the
Hodge bundle. In fact, $p^{*}({{\l}^{\otimes N}}_{\!\!m})$ and the pullback
${\pi}_T^{*}{{\tilde{\lambda}}_{m}}$ are also related by a canonical
isomorphism over $\Tg$.}


\noindent{\bf 3c. Functorial properties:}
We will require the following compatibility property for these
chosen isomorphisms.

Let
$$\matrix{Z&\mapright{f_1}&Y&\mapright{f_2} &X\cr}$$
be such that both $f_1$ and $f_2$ are covers,
of orders $n_1$ and $n_2$ respectively. Denote the genera of $Y$ and $Z$ by
$g_1$ and $g_2$ respectively. Let
$f_{1T}: {\cal T}_{g_1}\rightarrow\,{\cal T}_{g_2}$,
$f_{2T}:\Tg \rightarrow\,{\cal T}_{g_1}$ and
$(f_2\circ f_1)_T:\Tg \rightarrow\,{\cal T}_{g_2}$
be the induced maps of Teichm\"uller spaces.
If ${\l}_1$ and ${\l}_2$ denote the Hodge bundles on
${\cal T}_{g_1}$ and ${\cal T}_{g_2}$, respectively, and $\l$ the Hodge on
$\Tg$, then we have a {\it commutative diagram:}

\smallskip
$$
\matrix{{\l}^{n_1n_2}&\mapright{id}&{\l}^{n_1n_2}\cr
\mapdown{f_2^-}&&\mapdown{(f_2\circ f_1)^-}\cr
f_{2T}^*{\l}_{1}^{n_1}&\mapright{f_1^-}&(f_2\circ
f_1)_T^*{\l}_{2}\cr}
\leqno(3.4)
$$

In the diagram above the homomorphism of
Hodge bundles induced by any given covering has been denoted by using the
superscript $``-"$ over the notation for that covering.


\bigskip
\noindent {\underline{\it Section 4. HERMITIAN BUNDLE ISOMORPHISMS:}}

Recall the group of modular invariant line bundles on Teichm\"uller space,
$\Pth$, carrying hermitian structure, that was introduced in Section 1b.
We consider its ``rational version'' by tensoring (over $\Z$) with $\Q$:
$$
\Pthq:=  \Pth {\otimes}_{\Z}\Q
\leqno(4.1)
$$
It is to be noted that hermitian metrics, curvature forms (as members
of cohomology of the base with $\Q$-coefficients), etc. make perfect
sense for elements of rational Picard groups just as for usual
line bundles.

Lemma 1.1 now places things in perspective. We see from that result
that $\Pthq$ is canonically isomorphic to $\Pa \otimes \Q$ as well as
$\Ph \otimes \Q$.

We have seen in Section 1 that the Hodge bundle, and the
higher $DET$ bundles, admit their Quillen metric, such
that the curvature of the holomorphic Hermitian connection is a multiple of
the Weil-Petersson form. Therefore, these bundles $({\l}_m, Quillen)$
for $m=0,1,2,..$ constitute some interesting and canonical elements of
$\Pthq$.

As before, suppose that
$\pi:\,Y \rightarrow\,X$ is a covering that allows factorisation into
abelian covers, as in the set up of Theorem 3.2. Let the degree of the
covering be $N$.

First of all we note that the
isomorphism $\rho$ (Theorem 3.2) induces a modular invariant structure on
${\pi}_T^*{\tilde{\lambda}}$.  Note that though $\rho$
is defined only up to a non-zero constant, since multiplication by scalars
commutes with the modular invariance structure on $\l$, the modular
invariant structure obtained on ${\pi}_T^*{\tilde{\lambda}}$ does
not depend upon the exact isomorphism chosen.
{\it Thus we can put a
restriction on the isomorphism $\rho$, of Theorem 3.2, associated to
the covering $\pi$. We demand that it be
unitary with respect to the Quillen metric on
$({\psi})^*{\l}^N$ and the pull-back, by ${\pi}_k$,
of the Quillen metric on ${\tilde{\lambda}}$}.
In that case ${\bar F} = \rho$ is determined up to a scalar of norm one.

We have therefore proved the following Proposition that is fundamental
to our construction of Universal $DET$ bundles over inductive limits of
Teichm\"uller spaces:

\noindent {\bf Proposition 4.1:}
{\it Let $\pi$ be a covering with $N$ sheets as above. Then
$\l$ and ${\pi}^*_T((1/N){\tilde{\lambda}})$ coincide as elements of
$\Pth_{\Q}$. Here $(1/N){\tilde{\lambda}}$ is considered as an element of
$\Pttq$. By the same token, $\l_m$ and the pullback via ${\pi}_T$ of the
(1/N) times the $m^{th}$ $DET$ bundle over $\Tt$ coincide.}
\smallskip

\noindent {\it Remark}:~ In the diagram 3.4 (Section 3c), if all the
bundles are equipped with respective Quillen metric then the diagram becomes a
commutative diagram of Hermitian bundles. In other words, all the
morphisms in 3.4 are unitary.
\smallskip

More generally, let $(L,h)\in \Ptt$. Denote the pulled-back line bundle
by $L':= {\pi}_{T}^{*}L\in \Pt$. Because of modular invariance,
notice that $(L,h)$ becomes a hermitian line
bundle on $\Mt$. Consider the hermitian line bundle ${\pi}^*_k(L,h)$
on $\Mgk$, where ${\pi}_k$ is the induced map of Proposition 2.1.
The line-bundle ${\psi}^*L'$,
(where ${\psi}:\Mgk\rightarrow\,\Mg$ is the natural projection),
is canonically isomorphic
to ${\pi}^*_kL$. Now, using this isomorphism, and averaging the metric
${\pi}^*_kh$ over the fibers of $\psi$, a hermitian metric on the bundle
$L'$ over $\Mg$ is obtained. In other words $L'$ on $\Tg$ gets a
modular invariant hermitian structure. Thus we have a natural mapping:

$$
{\pi}_{T}^{\Q}~:\Pttq \longrightarrow \Pthq
\leqno(4.2)
$$

\smallskip
\noindent {\bf Proposition 4.2} {\it The natural homomorphism (4.2)
enjoys the following properties:

\noindent (1).
The image of ${\pi}_{T}^{\Q}$ on the Hodge bundle with Quillen
metric coincides with $N$ times the Hodge with Quillen metric over $\Mg$,
precisely as stated in Proposition 4.1. The
obvious corresponding assertion holds also for the higher DET bundles.

\noindent (2). For a pair of coverings
$$\matrix{Z&\mapright{f_1}&Y&\mapright{f_2} &X\cr}$$
as in the set up of 3.4, let  ${\pi}_{1}^{\Q}$
be the homomorphism given by (4.2) for the covering
$f_1$. Similarly, ${\pi}_{2}^{\Q}$ and ${\pi}^{\Q}$ be the
homomorphisms for $f_2$ and $f_2\circ f_1$ respectively. Then}
$${\pi}_{2}^{\Q}\circ {\pi}_{1}^{\Q}~ = ~ {\pi}^{\Q}$$


\bigskip
\noindent {\underline{\it Section 5. THE MAIN RESULTS:}}

\noindent{\bf 5a. Abelian lamination and its complex structures:}
Consider the inverse system of isomorphism classes of
finite covers of $X$,
which allow factorization into a sequence of {\it abelian} coverings.
The mappings in this directed system are those that commute
with the projections onto $X$.
One constructs thereby an interesting Riemann surface
lamination (``solenoidal Riemann surface'', see [S], [NS]),
by taking the inverse limit of this projective system
of surfaces.  This inverse limit, which we denote by
$H_{\infty, ab}$, is a compact space that is fibered over $X$
with a Cantor set as fiber. Each path component (``leaf'') of this
Riemann surface lamination is a covering space over the original
base surface $X$, the covering group being the commutator
subgroup of the fundamental group of $X$.

\noindent {\it Remark:} ~The general ``universal hyperbolic lamination''
of this type, where {\it all} finite coverings over $X$ are allowed,
has been studied in [NS]. That inverse limit space was denoted
$H_{\infty}$, and its Teichm\"uller space was shown to possess a
convergent Weil-Petersson pairing.

Let this projective system of covering mappings be
indexed by some set $I$, and for $i\in I$, let $p_i: X_i\rightarrow X$
be the covering and $Teich(p_i)=q_i: \Tg \rightarrow {\cal T}_{g_i}$
be the induced injective immersion between the Teichm\"uller spaces.
Thus, corresponding to the inverse system of surfaces we have a direct
system of Teichm\"uller spaces.

Let $g_i$ denote the genus of $X_i$. The direct limit space

$$\Tin\,:=\,\ \lim_{\rightarrow} {\cal T}_{g_i}
\leqno(5.1)
$$
of these Teichm\"uller spaces constitute the space of ``TLC''
[transversely locally constant] complex structures on the Riemann surface
lamination $H_{\infty, ab}^{n}$.

In fact, the lamination has a Teichm\"uller space parametrising
complex structures on it. This space is precisely the
completion in the Teichm\"uller metric  of the space $\Tin$. It is
the Teichm\"uller space of the abelian lamination, $T(H_{\infty, ab})$
and it is Bers-embeddable and possesses the  structure of an infinite
dimensional separable complex manifold. Notice that the various finite
dimensional Teichm\"uller spaces for Riemann surfaces of genus $g_i$
are all contained faithfully inside this
infinite dimensional ``universal'' Teichm\"uller space, from its very
construction.

\noindent
{\it Remark:}~We emphasize that the above construction is stable with
respect to change of the initial base surface $X$ . In fact, starting from
two distinct genera $X_1$ and $X_2$, the inverse limit solenoidal surface,
and the corresponding direct limit of Teichm\"uller spaces, will be
isomorphic whenever one finds a common covering surface $Y$ covering both
the bases by a (product of) abelian covers. That is easily done provided
the genera of $X_1$ and $X_2$ are both at least two. The families become
cofinal from $Y$ onwards. Hence the limiting objects are isomorphic.

We are interested in studying the space of complex
structures on this laminated surface, and constructing the canonical
$DET$ hermitian line bundle over this Teichm\"uller space,
{\it in order to obtain
the Mumford isomorphisms between the relevant $DET$ bundles over
this inductive limit space.}

\noindent
{\bf 5b. Line bundles on ind spaces:}
A line bundle on the inductive limit of an inductive system of varieties or
spaces, is, by definition ([Sha]), a collection of line bundles
on each stratum (i.e., each member of
the inductive system of spaces) together with compatible bundle maps.
The compatibility condition for the bundle maps is the obvious one relating
to their behaviour with respect to compositions, and guarantees that the
bundles over the strata themselves fit into an inductive system.
A bundle with hermitian metric is a collection with hermitian
metrics such that the connecting bundle maps are unitary.
Such a direct system of line bundles can clearly be thought of as
an element of the inverse limit of the Picard varieties of the
stratifying spaces. See [KNR] [Sha].

A ``rational''  hermitian line bundle over the inductive
limit is thus clearly an element of the inverse limit
of the $Pic \otimes \Q$ 's, namely of the space
$\lim_{\leftarrow}Pich({\cal T}_{g_i})_{\Q}$.  The following two
theorems were our chief aim.

\smallskip
\noindent
{\bf Theorem 5.1: Universal $DET$ line bundles:}~
{\it There exist canonical elements of the inverse limit
$\lim_{\leftarrow}Pich({\cal T}_{g_i})_{\Q}$, namely hermitian line bundles
on the ind space $\Tin$, representing the Hodge and higher DET bundles with
respective Quillen metrics:
$$
{\La}_{m} \in \lim_{\leftarrow}Pich({\cal T}_{g_i})_{\Q}, ~~m=0,1,2,..
$$
The pullback of $\La_m$ to each of the stratifying
Teichm\"uller spaces ${\cal T}_{g_i}$
is $(n_i)^{-1}$ times the corresponding Hodge or higher DET bundle
with Quillen metric over ${\cal T}_{g_i}$.}
\smallskip

\noindent
{\bf Proof:} The foundational work is already done in
Theorem 3.2 and Proposition 4.1 above. In fact,
let ${\l}_{0,i}$ represent the Hodge bundle with Quillen
metric in $Pich({\cal T}_{g_i})_{\Q}$. Then Proposition 4.1 implies
that, for $i\in I$ taking the element

$$(1/{n_i}){\l}_{0,i} \in Pich({\cal T}_{g_i})_{\Q}$$
provides us a compatible family of hermitian line bundles (in the rational
Pic) over the stratifying Teichm\"uller spaces -- as required in the
definition of line bundles over ind spaces.
The connecting family of bundle maps is determined (up to a scalar)
by Theorem 3.2.

The property that the connecting
unitary bundle maps for the above collection are {\it compatible}, follows
from Theorem 3.2 and the commutative diagram (3.4) of Section 3.
Notice that prescribing a base point in $\Tg$, and a vector of unit norm
over it, fixes uniquely all the scaling factor ambiguities in
the choices of the connecting bundle maps. We have therefore
constructed the universal Hodge, $\La_{0}$, over $\Tin$.

Naturally, the above analysis can be repeated verbatim for
the higher d-bar families, and one thus obtains elements
$\La_{m}$ for each positive integer $m$. Again the pullback of
$\La_{m}$ to any of the stratifying ${\cal T}_{g_i}$ produces
$(n_i)^{-1}$ times the $m$-th $DET$ bundle with Quillen metric
living over that space.
$\hfill{\Box}$

\smallskip
\noindent
{\bf Theorem 5.2: Universal Mumford isomorphisms:}~
{\it Over the direct limit Teichm\"uller space $\Tin$ we have
$$\La_{m}\,=\,(6m^{2} + 6m +1){\La_{0}}$$ as an equality of
hermitian line bundles.}
\smallskip

\noindent
{\bf Proof:} Follows directly from the genus-by-genus isomorphisms
of (1.11), and our universal line bundle construction of this paper.
$\hfill{\Box}$

\noindent
{\bf 5c. Polyakov measure on $\Mg$ and our bundles:}
The quantum theory of the closed bosonic string is the theory of a sum
over random surfaces (``world-sheets") swept out by strings propagating
in Euclidean spacetime $\R^{d}$.  Owing to the conformal invariance
property enjoyed by the Polyakov action,
that summation finally reduces to integrating out a certain measure -
called the Polyakov measure - on the moduli space of conformally
distinct surfaces, namely on the {\it{parameter spaces of Riemann
surfaces}}.  It is well-known ([P], [Alv], [BM]) that the original
sum over the infinite-dimensional space of random world-sheets (of
genus g), reduces to an integral over the finite-dimensional  moduli
space, ${\cal{M}}_g$, of Riemann surfaces of genus $g$ precisely when
the background spacetime has dimension $d = 26$.
In the physics literature, this is described by saying that the
``conformal anomaly" for the path-integral vanishes when $d=26$.

That magic dimension 26, at which the conformal anomaly cancels, can also
be explained as the dimension in which the ``holomorphic anomaly" cancels.
The vanishing of the holomorphic anomaly at $d = 26$ can be interpreted as
the statement of Mumford's isomorphism (1.11) above for the case $m=1$.
In fact, that result shows, in combination with out interpretation
of the $DET$ bundles $\l$ and ${\l}_{1}$ (see Section 1c),
that the holomorphic line bundle $K \otimes {\l}^{-m}$ over
${\cal{M}}_g$ is the trivial bundle when $m = d/2=13$.
Here $K={\l}_1$ denotes the canonical line bundle, and $\l$ the
Hodge bundle, over ${\cal{M}}_g$.

The Polyakov volume form on ${\cal{M}}_g$ can now be given the
following simple interpretation. Fixing a volume form (up to
scale) on a space amounts to fixing a fiber metric (up to scale) on the
canonical line bundle over that space.  But the Hodge bundle $\l$
has its natural Hodge metric (arising from the $L^2$ pairing of
holomorphic 1-forms on Riemann surfaces).  Therefore we may transport
the corresponding metric on ${\l}^{13}$ to $K$ by Mumford's isomorphism,
(as we know the choice of this isomorphism is unique up to scalar) --
thereby obtaining a volume form on ${\cal{M}}_g$.
[BK] showed that this is none other than the Polyakov volume.
Therefore, the presence of Mumford isomorphisms over the
moduli space of genus $g$ Riemann surfaces
describes the Polyakov measure structure thereon.

Above we have succeeded in fitting together the Hodge and higher
$DET$ bundles over the ind space $\Tin$, together with the relating
Mumford isomorphisms. We thus have from our results a structure
on $\Tin$ that qualifies as a genus-independent, universal, version of
the Polyakov structure.

\noindent{\it{Remark}}:~ Since the genus is considered the
perturbation parameter in the above formulation of the standard
perturbative bosonic Polyakov string theory, we may consider our work as
a contribution towards a {\it non-perturbative} formulation of that
theory.

\noindent{\it{Remark}}:~ In [NS] we had shown the existence of a convergent
Weil-Petersson pairing for the Teichm\"uller space of the universal
lamination $H_{\infty}$. Now, Polyakov volume can be written as
a multiple by a certain combination of Selberg zeta functions
of the Weil-Petersson volume. Using therefore
Theorems 5.1 and 5.2 in conjunction with the
Weil-Petersson on $T(H_{\infty,ab})$, we have an interpretation of that
combination of Selberg zetas as an object fitted together
over each stratifying Teichm\"uller space.

It should be noted in this context that the line bundles over $\Tin$ that
we have constructed do not appear as restrictions to the strata
from any rational line bundle over the completed space $T(H_{\infty,ab})$.

\bigskip
\noindent
{\underline{\it Section 6: EXTENSION TO MODULI WITH MARKED POINTS:}}

\noindent
{\bf 6a. Moduli spaces with $r$ distinguished points:}
As always, denote by $X$
a compact oriented $C^{\infty}$ surface of genus $g$, and
$$S=\{p_1,\dots ,p_r\}$$ a finite ordered subset of $r$ distinct
points on $X$.

Let ${\cal M}_g^r$ stands for the moduli space of $r$-pointed compact
Riemann surfaces of genus $g$. By definition, it is the space of
$Diff^{+}(X,\{p_1,\dots ,p_r\})$-orbits of conformal
structures on $X$. The corresponding Teichm\"uller space
is denoted by $\Tgr$.

Fix any integer $k \geq 2$. Let $\Mgrk$ be the moduli space with
level $k$. In other words, $\Mgrk$ is the moduli of triplets of the form
$(M,\Ga,\r)$, where $M$ is a Riemann surface of genus $g$, $\Ga$ is
an ordered subset of cardinality $r$ and $\rho$ is a basis of the
$\Z/k\Z$ module $H_1(M-\Ga ,\Z/k\Z)$. $\Mgrk$ is a complex manifold
of dimension $3g-3+r$ covering the orbifold $\Mgr$.

\noindent {\it Generators for Pic:}
Of course, the Hodge line bundle can be introduced as before as an
element of $\Pfr$. Moreover, for each $i\in \{1,\ldots, r\}$, there
is a line bundle $L_i$ on $\Mgr$, whose fiber
$L_{iX}$, over $X\in\Mgr$, is $K_{p_i}$ (recall $p_i$ is the $i$-th
element of $S$). Here $K_x$ denotes the fiber of the canonical bundle
at $x$ on any Riemann surface.
The Picard group of of the moduli functor for
$\Mgr$, $\Pfr$, is known to be generated by the Hodge
bundle together with these $\{L_i\}$. See [AC].

\noindent {\bf 6b. Construction of morphisms between moduli:}
Let $$\pi: \X \longrightarrow X$$ be a Galois cover, with cyclic
Galois group being $\Z/l\Z$, where $l$ is a {\it prime} number. {\it We
assume that $\pi$ is ramified exactly over $S$}. Since the Galois cover
is of prime order, it is easy to see that the map $\pi$ is
necessarily totally ramified over each point of $S$. Let the
genus of $\X$ be $\tilde g$.

Using the same techniques as for Proposition 2.1, we can show the following:

\smallskip
\noindent {\bf Proposition 6.1:}~{\it Let the level number $k$ be
fixed at any multiple of the number of sheets $l$. Then the covering
$\pi$ induces natural maps

$$(i){\pi}_T: \Tgr\,\longrightarrow\,{\Tt^{r}}$$

$$(ii){\pi}_k: \Mgrk\,\longrightarrow\,{\Mt}^{r}$$
Moreover if
$q:\Mgrk\longrightarrow\,{\cal M}_g^r[l]$ be the projection induced by
the homomorphism $\Z/k\rightarrow \Z/l$ then the diagram
corresponding to the one shewn in Proposition 2.1 again commutes.}
\smallskip

Now let $\pi: Z \longrightarrow\,X$ be an {\it abelian} Galois covering
which is totally ramified over the subset $S$, and nowhere else.
Exactly the same results as proved in Sections 3 and 4 above
go through with essentially no extra trouble. The set of points
of ramification at any covering stage remains of cardinality $r$, and
one needs to keep track of this subset all along.

Indeed, if we identify, as in Section 4, the hermitian
bundles on $\Mgr$ with modular invariant bundles on $\Tgr$, then
associated to the covering $\pi$ there is
a morphism which is precisely the analog of morphism (4.2):
$$
{\pi}_{T}^{r,\Q}: {\Pttqr} \longrightarrow {\Pthqr}
\leqno(6.1)
$$
(6.1) exists because, as we have seen for the Hodge bundles,
the bundles $L_i$ are also
well-behaved under pull-back by $\pi$ (again the $N$-th power
of $L_i$ identifies with ${\tilde L}_{i}$). As we mentioned above,
Hodge together with
these $r$ bundles generate $\Pfr$. Therefore (6.1) is defined and
well-behaved. This homomorphism satisfies the analogue of Proposition 4.2.

\noindent
{\bf 6c. Line bundles on the inductive limit $\Tir$:}
Consider the inverse system of isomorphism classes of
finite covers of $X$, totally ramified over the finite subset $S$,
which allow factorization into a sequence of {\it abelian} coverings.
The mappings in this directed system are those that commute
with the projections onto $X$. As in the case without punctures,
one constructs thereby an interesting Riemann surface
lamination by taking the inverse limit of this projective system
of surfaces.  This inverse limit, which we denote by
$H_{\infty, ab}^{r}$,
is a compact space that is fibered over $X$ with a Cantor set as fiber.

Let this projective system of covering mappings be
indexed by the set $I$, and for $i\in I$, let $p_i: X_i\rightarrow X$
be the covering and let $q_i: \Tgr \rightarrow {\cal T}_{g_i}^r$
be the induced map between Teichm\"uller spaces. Here $g_i$
denotes the genus of $X_i$ and $n_i$ is the order of the covering.
As before, this constitutes an inductive
system of Teichm\"uller spaces, and we have the direct limit

$${\Tir} := \lim_{\rightarrow} {\cal T}_{g_i}^r$$

\noindent
of these Teichm\"uller spaces constituting the space of ``TLC''
[transversely locally constant] complex structures on the Riemann surface
lamination $H_{\infty, ab}^{r}$.

The completion in the Teichm\"uller metric of the space $\Tir$ is
the Teichm\"uller space of the $r$-pointed abelian lamination;
it is denoted $T(H_{\infty, ab}^{r})$.
$\Tir$ itself comprises the space of ``TLC''
complex structures on the Riemann surface
lamination $H_{\infty, ab}^{r}$. Notice that the various finite
dimensional Teichm\"uller spaces for Riemann surfaces of genus $g_i$ with
$r$ distinguished points are all contained faithfully inside this
infinite dimensional ``universal'' Teichm\"uller space.

The main result, Theorem 5.1 above, as well as its proof, goes through
and has the corresponding statement for each of the $DET$ bundles.
It should be noted that the Pic-generator bundles $L_i$, $i=1,..,r$,
also fit together coherently over the ind space $\Tir$ (as members
of its rational Pic).

\smallskip
\noindent
{\bf Theorem 6.2: Universal $DET$ line bundles:}~
{\it There exist canonical elements of the inverse limit
$\lim_{\leftarrow}Pich({\cal T}_{g_i}^{r})_{\Q}$, namely hermitian line bundles
on the ind space $\Tir$, representing the Hodge and higher DET bundles with
respective Quillen metrics.

The pullback of $\La_m$ to each of the stratifying
Teichm\"uller spaces ${\cal T}_{g_i}^{r}$
is $(n_i)^{-1}$ times the corresponding Hodge or higher DET bundle
with Quillen metric.}

\newpage

\begin{center}
{\bf Appendix: Discreteness of the Picard variety of moduli spaces with
marked points and level structure}
\end{center}
\bigskip

\noindent
Professor Eduard Looijenga (Utrecht) has sent us the following argument for:

\noindent{\bf Theorem A:} ~{\it The Picard group
of the moduli space of pointed Riemann surfaces with a principal level
structure is discrete. In fact, the stronger assertion that the corresponding
mapping class group has trivial first cohomology (with $\Z$ coefficients)
holds.}

\noindent {\it Remark:} ~As mentioned before, see also [Hain].
The above result is
clearly of independent interest and reproves our basic Theorem 2.3, namely
that ``curvature determines the bundle'' on such moduli spaces.

\noindent {\bf Proof:}
Let us be given a compact connected oriented surface $S$ of
genus $g$. Write $V_g$ for $H^1(S)$ and denote by
$\Gamma _g$ resp. $Sp_g$ the mapping class group of $S$  resp. the
group of integral symplectic transformations of $V_g$. There is a natural
surjection $\Gamma _g\to Sp_g$ whose kernel $T_g$ is known as the Torelli group
of $S$. In a series of papers Dennis
Johnson proved that if $g\ge 3$ (an assumption in force from now on), $T_g$ is
finitely generated and that $H^1(T_g)=Hom (T_g,\Z)$ can be naturally identified
with the quotient of $\wedge ^3 V_g$ by $\omega\wedge V_g$, where
$\omega\in \wedge ^2 V_g$ corresponds to the intersection product.
(This is precisely the primitive cohomology in degree 3
of the Jacobian $H_1(S;{\R}/{\Z})$ of $S$, but that fact
will not play a role here.) This isomorphism is $Sp_g$-equivariant and so there
are no nonzero elements in $H^1(T_g)$ invariant under $Sp_g$. This is also
true for any subgroup of $Sp_g$ of finite index; in particular it is so for the
principal level $n$ congruence subgroup $Sp_g[n]$. If we define the subgroup
$\Gamma _g[n]$ of $\Gamma _g$ by the exact sequence
$$
0\to T_g\to \Gamma _g[n]\to Sp_g[n]\to 0,
$$
then what we want is the vanishing of $H^1(\Gamma _g[n])$. The spectral
sequence associated to this exact sequence
$$
E^{p,q}_2=H^p(Sp_g[n];H^q(T_g))\Rightarrow H^{p+q}(\Gamma _g[n])
$$
gives in degree 1 the short exact sequence
$$
0\to H^1(Sp_g[n])\to H^1(\Gamma _g[n])\to H^1(T_g)^{Sp_g[n]}\to
$$
We noted already that the term on the right is trivial. So it remains to see
that $H^1(Sp_g[n]) \cong H^1(Sp(2g,{\Z})[n])$ is trivial. A. Borel has
shown that in this range the real cohomology of
$Sp(2g,{\Z})[n]$ can be represented by
translation-invariant forms on the corresponding Siegel space. It is therefore
trivial, too. (The result we are alluding to is in Theorem 7.5 in his
paper [Bo], although sharper results have been obtained by him since.The
constants involved in applying Theorem 7.5 can be noted from [KN] (reference
[16]
of [Bo]). In the case at hand, both constants ``$c$'' and ``$m$''that
are involved are at least $1$, and hence our argument is valid,
provided $g$ is at least $3$. )

\smallskip
In order to obtain the corresponding result for the moduli space of $k$-pointed
genus $g$ curves we fix some notation: let $z_1,z_2,z_3,\dots$ be a
sequence of distinct points of $S$ and put $V_g^k:=H^1(S-\{z_1,\dots ,z_k\})$.
Notice that the restriction maps give inclusions $V_g=V_g^0=V^1_g\subset
V_g^2\subset\cdots$. Denote the group of automorphisms of $V_g^k$
that leave $V_g\subset V_g^k$ invariant, act symplectically on $V_g$ and are
the
identity on $V_g^k/V_g$ by $Sp_g^k$ and let
$\Gamma _g^k$ stand for the connected
component group of the group $Diff ^+(S,\{z_1,\dots ,z_k\})$ of orientation
preserving diffeomorphisms of $S$ that leave
$z_1,\dots ,z_k$ pointwise fixed. It
is known that the natural homomorphism $\Gamma _g^k\to Sp_g^k$ is surjective.

We prove with induction on $k$ that $H^1(\Gamma _g^k[n])$ is trivial.
So assume $k\ge 1$ and the claim proved for $k-1$.
We have an exact sequence
$$
1\to \pi _1(S-\{z_1,\dots ,z_{k-1}\},z_{k})\to \Gamma _g^k[n]\to \Gamma
_g^{k-1}[n]\to 1. $$
Its associated spectral sequence for cohomology gives the short exact sequence
$$
0\to H^1(\Gamma _g^{k-1}[n])\to H^1(\Gamma _g^k[n])\to
(V_g^{k-1})^{Sp_g^{k-1}}\to
$$
The left-hand term is trivial by assumption and it is easy to see that the
right-hand term is trivial, too. This gives the vanishing of the middle term.

Let ${\cal M}_g^k$ stands for the moduli space of $k$-pointed compact
Riemann surfaces of genus $g$ (Section 6a).
and let ${\cal M}_g^k[n]$ mean the moduli with level $n$. We see from above
that for $g\ge 3$, ${\cal M}_g^k[n]$ has vanishing first
Betti number and has therefore a discrete Picard group, as desired.
$\hfill{\Box}$

\bigskip

\begin{center}
{\bf References}
\end{center}

\begin{enumerate}

\item[{[Alv]}] O. Alvarez, Theory of strings with boundaries:
fluctuations, topology and quantum geometry, {\it Nucl. Phys, B216}, (1983),
125-184.

\item[{[AC]}] E. Arbarello, M. Cornalba, The Picard group of the
moduli spaces of curves. {\it Topology, 26}, (1987), 153-171.

\item[{[BM]}] A.A. Beilinson, Yu. I. Manin, The Mumford form and the
Polyakov measure in string theory, {\it Comm. Math. Phys., 107}, (1986),
359-376.

\item[{[BK]}] A. Belavin, V. Knizhnik, Complex geometry and quantum
string theory. {\it Phys. Lett., 168 B}, (1986), 201-206.

\item[{[BF]}] J. Bismut, D. Freed, The analysis of elliptic families
:metrics and connections on determinant bundles. {\it Comm. Math. Phys.,
106}, (1986), 159-176.

\item[{[Bo]}]: A. Borel, Stable real cohomology of arithmetic groups,
{\it Ann. Sci. Ecole. Norm. Sup., 7}, (1974), 235-272.

\item[{[Bos]}] J.B. Bost, Fibres determinants, determinants regularises
et mesures sur les espaces de modules des courbes complexes,
{\it Semin. Bourbaki, 152-153}, (1987), 113-149.

\item[{[G]}]: T. Gendrone, Thesis, CUNY Graduate Center, (to appear).

\item[{[Hain]}] R. Hain,  Torelli groups and the geometry of moduli
spaces of curves, preprint, (1993-94).

\item[{[Har]}] J. Harer, The second homology of the mapping class group
of an orientable surface, {\it Invent. Math., 72}, (1983), 221-239.

\item[{[H]}] R. Hartshorne, {\it Algebraic Geometry}, Springer Verlag,
(1977).

\item[{[HM]}] J.Harris and D. Mumford, On the Kodaira dimension of the
moduli space of curves, {\it Invent. Math., 67}, (1982), 23-86.

\item[{[KN]}] S. Kaneyuki and T. Nagano, Quadratic forms related to
symmetric spaces, {\it Osaka Math J., 14}, (1962), 241-252.

\item[{[KNR]}] S. Kumar, M.S. Narasimhan, A. Ramanathan, Infinite
Grassmannians and moduli spaces of $G$-bundles, {\it Math.
Annalen}, (to appear).

\item[{[M]}] C.C. Moore, Group extensions of p-adic and adelic
linear groups.{\it Publ. Math. I.H.E.S.,  35}, (1968), 5-70.

\item[{[Mum]}] D. Mumford, Stability of projective varieties, {\it Enseign.
Math., 23}, (1977), 39-100.

\item[{[NS]}] S. Nag, D. Sullivan, Teichm\"uller theory and the universal
period mapping via quantum calculus and the $H^{1/2}$ space on the
circle, {\it Osaka J. Math., 32, no.1}, (1995), (to appear).
[Preprint prepared at I.H.E.S., 1992.]

\item[{[P]}] A.M. Polyakov, Quantum geometry of bosonic strings, {\it Phys.
Lett., 103B}, (1981), 207-210.

\item [{[Q]}] D. Quillen, Determinants of Cauchy-Riemann operators
over a Riemann surface. {\it Func. Anal. Appl., 19}, (1985), 31-34.

\item[{[Sha]}] I.R. Shafarevich, On some infinite-dimensional groups II,
{\it Math USSR Izvest., 18}, (1982), 185-194.

\item[{[S]}] D. Sullivan, Relating the universalities of Milnor-Thurston,
Feigenbaum and Ahlfors-Bers, in Milnor Festschrift, {\it ``Topological Methods
in Modern Mathematics''}, (eds. L.Goldberg and A. Phillips) Publish or
Perish, (1993), 543-563.

\item[{[ZT]}] P.G. Zograf and L.A. Takhtadzhyan, A local index theorem for
families of {$\bar \partial$}- operators on Riemann surfaces,
{\it Russian Math Surveys, 42}, (1987), 169-190.

\end{enumerate}

\end{document}